\documentclass[aps,letterpaper,preprint,amssymb,pra]{revtex4} 
\usepackage[none]{hyphenat} 
\usepackage{bm}
\usepackage{amsmath}
\usepackage{amssymb}
\usepackage{times}
\usepackage{latexsym}
\usepackage{graphicx}
\usepackage{color}
\usepackage{subfigure}
\begin{document}

\title{ Kaleidoscopic Symmetries and Self-Similarity of Integral Apollonian Gaskets}
\author{  Indubala I Satija }
\affiliation{Department of Physics, George Mason University , Fairfax, VA 22030, USA}
\date{\today}
\begin{abstract}

{ We describe various kaleidoscopic and  self-similar aspects of the integral Apollonian gaskets -   fractals consisting of close packing of circles with integer curvatures.  Self-similar recursive structure of the whole gasket is shown to be encoded in transformations that forms the modular group $SL(2,Z)$.
The asymptotic scalings of curvatures of the circles are  given by a special set  of quadratic irrationals with continued fraction   $[n+1: \overline{1,n}]$ - that is  a  set of irrationals with period-2 continued fraction consisting of $1$ and  another integer $n$. Belonging to the class $n=2$, there exists a nested set of self-similar kaleidoscopic patterns that exhibit three-fold symmetry. Furthermore, the even $n$ hierarchy is found to mimic the recursive structure of the  tree  that generates all  Pythagorean triplets }
\end{abstract}

\maketitle

\begin{figure}[htbp] 
 \includegraphics[width = .6\linewidth,height=0.55 \linewidth]{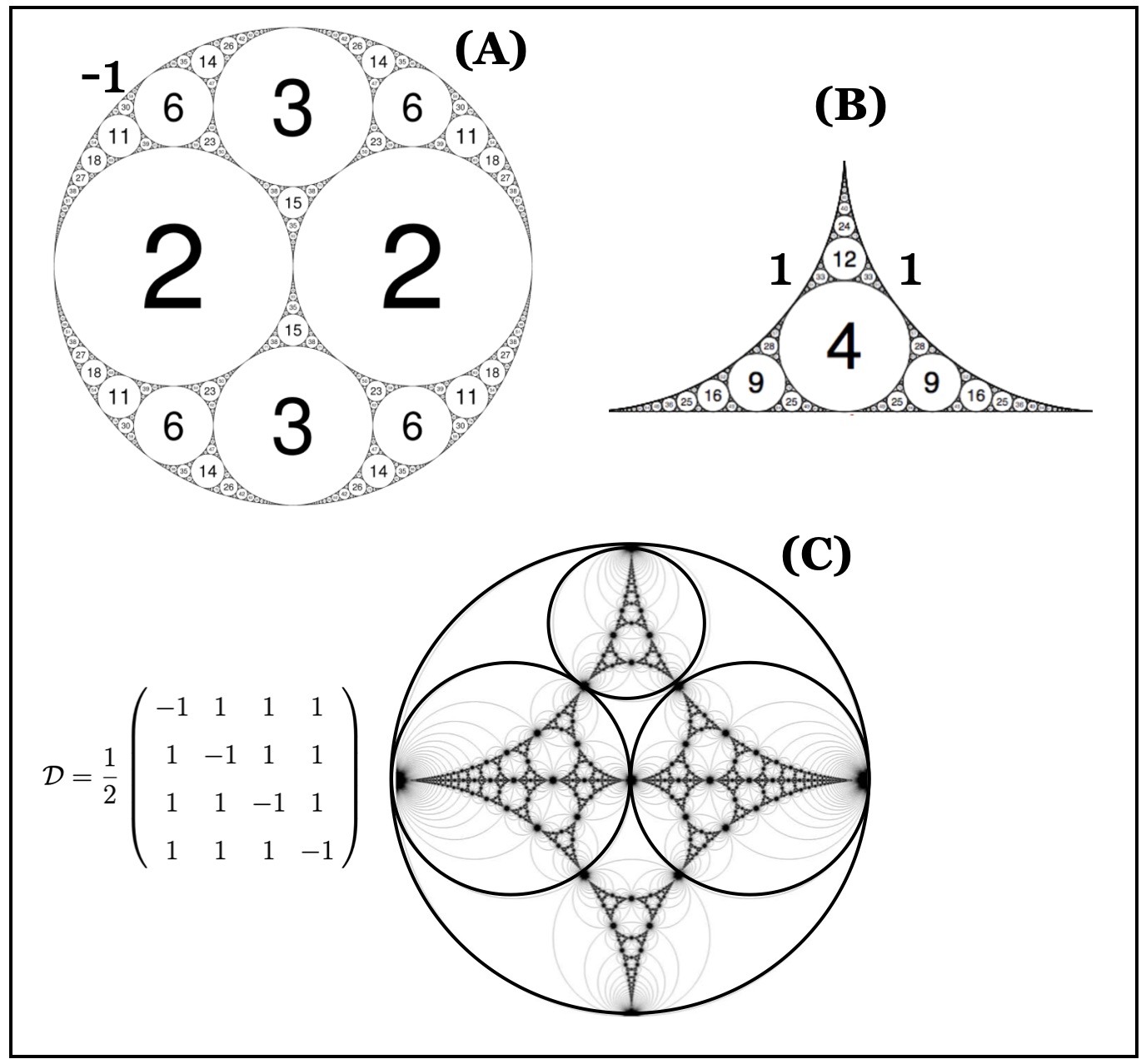}
\leavevmode \caption{ Integral Apollonian gaskets. The  gaskets shown in panels (A) and (B) are dual to each other,
as illustrated in panel (C) where we overlay the two gaskets. The duality described by the matrix $\mathcal{D}$ ( displayed in lower left ) corresponds to each circle in the dual set passing through three of the tangency points of the original set of circles, and the reverse holds as well. For example, $(-1,2,2,3)^T= \mathcal{D} ( 4, 1, 1, 0)^T$ and $(4,1,1,0) = \mathcal{D} ( -1, 2, 2, 3)^T$.}
\label{iagall}
\end{figure}

Integral Apollonian gaskets$(\mathcal{IAG})$\cite{IAG}  such as those shown in figure ~(\ref{iagall}) consist of close packing of circles of integer curvatures (reciprocal of the radii),
 where every circle is tangent to three others.
These are fractals where the whole gasket is like a kaleidoscope reflected again and again through an infinite collection of curved mirrors\ that  encodes  fascinating geometrical and number theoretical concepts\cite{AP}. The central themes of this paper are the kaleidoscopic and self-similar recursive properties described within the framework of M\"{o}bius transformations that maps circles to circles\cite{Cmap}.

Named in honor of Apollonius of Perga who studied the geometrical problem of mutually tangent circles before 300BC, 
the building block of Apollonian gaskets are configurations of four mutually tangent circles whose curvatures  $(\kappa_1, \kappa_2, \kappa_3, \kappa_4)$ satisfy the following relation\cite{prime, J1}:

\begin{equation}
 2( \kappa_1^2+\kappa_2^2+\kappa_3^2+\kappa_4^2) =  ( \kappa_1+\kappa_2+\kappa_3+\kappa_4)^2 .
 \label{Q4}
 \end{equation}
 Here, the curvature of the outer circle that encloses all inner circles has to be taken negative to satisfy an equation. 
Eq. (\ref{Q4})), known as  the  Descartes theorem, can also be written as
$ v \mathcal{D} v^T=0$.
 Here $ v = ( \kappa_1, \kappa_2, \kappa_3, \kappa_4)$ and $v^T$ is its transpose and the matrix $\mathcal{D}$ is defined in Fig. ~(\ref{iagall}).
 Configuration of circles satisfying Descartes theorem will be referred as the ``Descartes configurations".
 Originally discussed by French philosopher and mathematician René Descartes  in $1643$ in a letter to Princess Elizabeth of Bohemia,
  Apollonian circles were rediscovered in $1936$ by Chemistry Nobel laureate  Frederick Soddy who published in Nature a poetic version of Descartes’ theorem, which he called ``The Kiss Precise”.
    
 For a given integral Apollonian packing, the investigations about its Diophantine properties\cite{J2} such as what integers appear as curvatures and the  number of circles with prime curvatures have piqued many mathematicians\cite{prime}. The formula for the number of primes in an integral Apollonian gasket bears a striking resemblance to the Prime Number Theorem.\cite{prime}. These questions have been intimately linked to the hidden symmetries and consequently the deceptively simple geometrical construction  of close packing of circles is strongly tied with sophisticated mathematics involving group theory.  
 Furthermore, recent studies have
shown that these abstract fractals are related to a quantum fractal known as the Hofstadter butterfly \cite{Hof} which models  all possible integer quantum Hall states\cite{QHE} -  the exotic topological states of matter\cite{book, IIS2, IIS1}. A detailed investigation of the recursive nature of the  integral Apollonian gasket using M\"{o}bius transformations that form modular group $SL(2,Z)$ as described here 
shows that in close analogy to the integer curvatures of the gasket,
the recursive structure of the topological quantum numbers of the butterfly, explained earlier\cite{SW}  using quantum mechanics are in fact rooted in pure geometry and number theory.

\begin{figure}[htbp] 
 \includegraphics[width = .63 \linewidth,height=0.3 \linewidth]{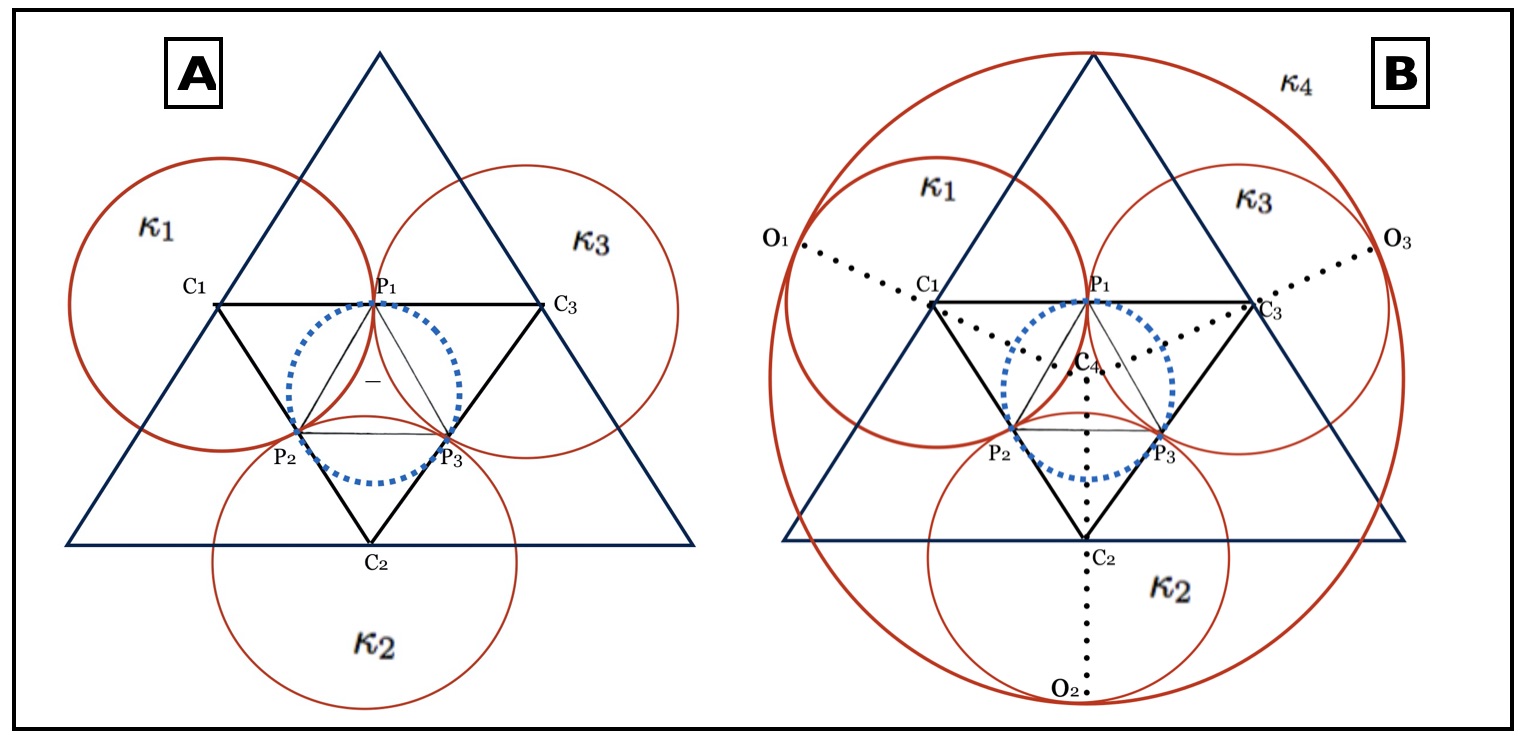} 
\leavevmode \caption{ (A) Given any three points $(P_1, P_2, P_3)$, one can construct  a configuration of three mutually tangent circles of curvatures $(\kappa_1, \kappa_2, \kappa_3)$: we first draw a circle ( dotted blue ) through these points and then draw tangents to this circle. Three tangents  meet at $C_1$, $C_2$ and $C_3$ and using $( C_1, C_2, C_3)$ as the centers, we now draw three circles which are mutually tangent as the distance between
$P_1$ and $C_1$ is same as the distance between $P_2$ and $C_2$ and so on.  The  curvature $\kappa_4$ and hence the radius of the fourth circle is determined by Descartes theorem ( Eq. (\ref{Q4})) . The center $C_4$ is determined in terms of the curvatures of the four circles\cite{4circle}.  }
\label{P3}
\end{figure}

 An exercise in geometry as illustrated in Figure (\ref{P3}), shows  how to obtain a configuration of four mutually tangent circles starting with any three points. As every Descartes configuration is fully determined by three distinct points in a plane,  this implies that any two such configurations are related by a unique conformal map - a M\"{o}bius transformation\cite{Cmap} of the complex plane of the form,
  \begin{equation}
   {\displaystyle z\mapsto  f(z)={\frac {az+b}{cz+d}}}
   \label{MM}
   \end{equation}
  A brief review of some of the properties of these transformations is given in Appendix A.   Section 1  illustrates kaleidoscopic properties of the $\mathcal{IAG}$.
In Section 2,   we introduce the Ford Apollonian gasket in which each circle’s curvature is the square of some non-negative integer.
 Its recursive structure is described in Section 3,
   by a M\"{o}bius transformation that forms $SL(2,Z)$ group where $(a,b,c,d)$ in the map (\ref{MM}) are integers and Section (4) discusses its scaling  properties.
   In section 5, M\"{o}bius maps that belongs to the group
  $SL(2,C)$ where  $(a,b,c,d)$  are in general complex numbers, are shown to map the entire $\mathcal{IAG}$ to the Ford Apollonian gasket. Section 6 shows one of a kind self-similar kaleidoscopic patterns of the gasket. In Section 7,  the hierarchical pattern of the Pythagorean tree is shown to encode part of the recursions of the $\mathcal{IAG}$. In section (8) we discuss the Apollonian group and its relation to modular group $SL(2,Z)$. Appendix B gives a proof of Descartes theorem for Ford Apollonian configurations where it is shown to be a consequence of a  property of the homogeneous functions. 

 \section{ Apollonian Gasket - A Kaleidoscope}
Apollonian gaskets exhibit numerous  kaleidoscopic symmetries - a rare and fascinating feature among known functions. In various figures below, the circular mirrors  describing kaleidoscopic symmetries
 will be shown as dotted blue circles.

 The existence of a simplest  kaleidoscope in a gasket follows from the Descartes formula (\ref{Q4}).
Given any three mutually tangent circles of curvatures $(\kappa_1, \kappa_2, \kappa_3)$, there are exactly two possible circles $\kappa_4$ or $\bar{\kappa}_4$  that are tangent to these  three circles:
\begin{equation}
\kappa_{4} ( \bar{\kappa}_4) =(\kappa_1+\kappa_2+\kappa_3) \pm 2 \sqrt{ \kappa_1 \kappa_2+\kappa_2 \kappa_3+ \kappa_3 \kappa_2}
\label{m4}
\end{equation}

Equation (\ref{m4}) leads to a linear equation connecting the two solutions:
  \begin{equation}
 \kappa_4 + \bar{\kappa}_4 = 2 ( \kappa_1+\kappa_2+\kappa_3) .
 \label{l4}
 \end{equation}
 
 This linear equation implies that if the original four circles have integer curvature, all of the circles in the packing will have integer curvatures. This is because  starting with three mutually tangent circles, we can  construct two distinct quadruplets $(\kappa_1, \kappa_2,\kappa_3, \kappa_4) $ and $(\kappa_1, \kappa_2,\kappa_3,  \bar{\kappa}_4)$ and thus adding additional circles to the gasket.
 In fact we obtain the entire gasket as the equation (\ref{l4}) can also be written as $ \kappa_2 + \bar{\kappa}_2 = 2 ( \kappa_1+\kappa_4+\kappa_3)$ and so on. 
 
   \begin{figure}[htbp] 
 \includegraphics[width = .6 \linewidth,height=0.57 \linewidth]{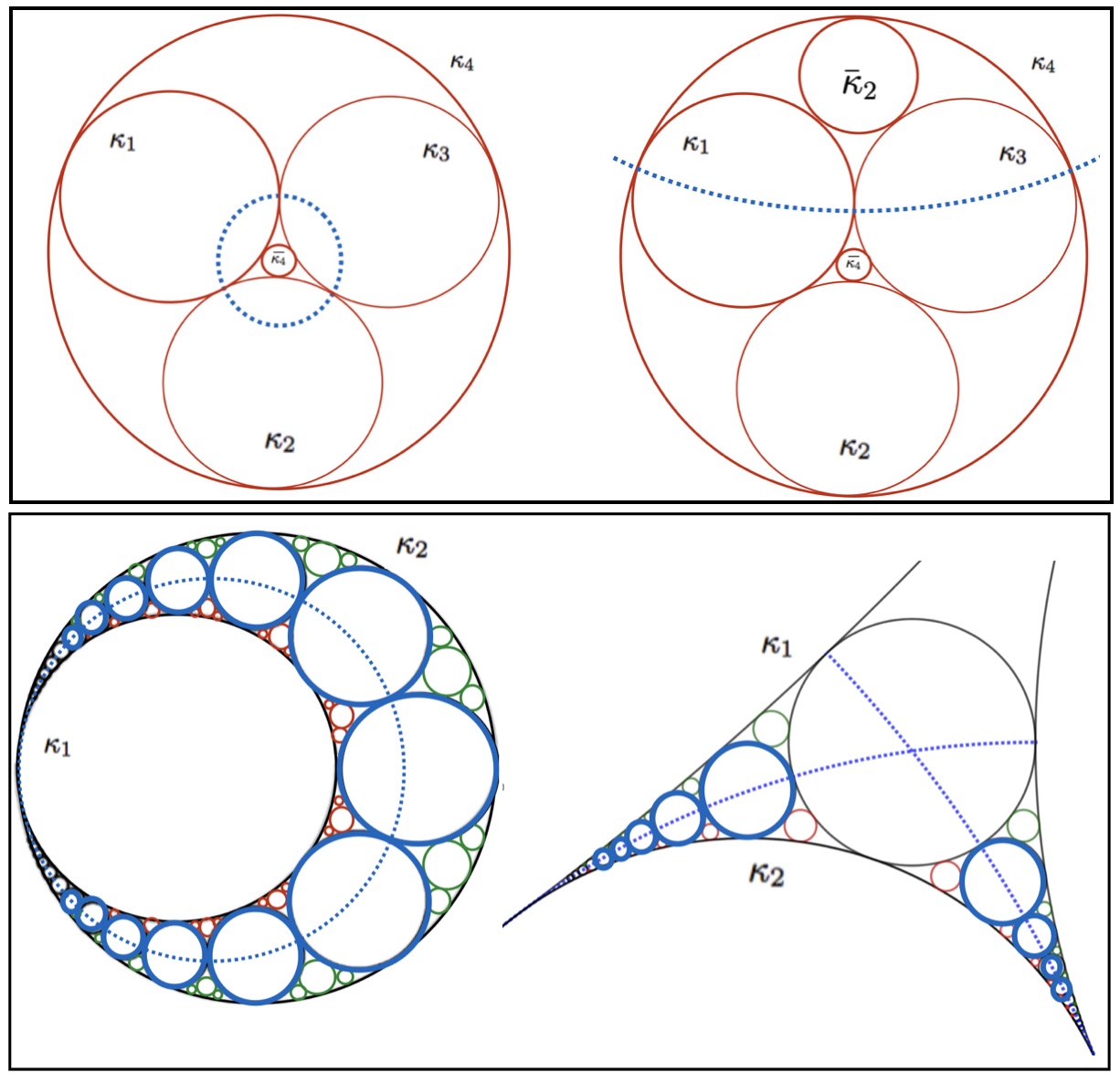} \leavevmode \caption{ Examples of circular mirrors shown as dotted blue lines. Upper panel shows that the two solutions of (\ref{l4}) are mirror image. The lower two panels show two examples: Pappus chain - chain of (blue) circles that are tangent to two mutually tangent (black) circles $(\kappa_1,\kappa_2)$ which are mirror images of each other through Pappus mirror ( ( dotted blue circles) . 
 Pappus mirror  passes through the tangency points of Pappus chain and 
 it reflects $\kappa_1$ and all the circles that are tangent to it ( red circles) to $\kappa_2$ and all the circles that are tangent to it ( green circles ).}
\label{PM}
\end{figure}

 A remarkable aspect of the pair of solutions such as $(\kappa_4 , \bar{\kappa}_4)$  is that they are mirror images of each other 
 through a circular mirror that passes through the tangency points of $(\kappa_1, \kappa_2, \kappa_3)$  as shown in Fig. (\ref{PM}).
 Therefore, packing of circles in the whole gasket is via kaleidoscopic images in the curvilinear triangular spaces in between the circles. 

Lower caption in Fig. ~(\ref{PM}) shows a different kind of mirror due to  global symmetries where an entire hierarchical pattern and  its mirror image is part of the gasket.  Given two
tangent circles  of curvatures $\kappa_1$ and $\kappa_2$, there exists a chain of (blue) circles where two consecutive circles of the chain are tangent to each other  as seen in the figure. These 
are examples of {\it Pappus chains}, investigated by Pappus of Alexandria in the 3rd century \cite{Pchain}. The tangency points of the chains of circles lie on a circle. The circular mirror at these tangency points - dubbed {\it Pappus mirror}"  reflect $\kappa_1$  into $\kappa_2$ and also the entire hierarchical set of circles that are tangent to $\kappa_1$  ( red circles) to the 
hierarchical set that are tangent to $\kappa_2$ ( green circles). In other words, for the hierarchical configuration, where all circles share a common circle, which we refer as the ``boundary circle", there exists twin configuration - its mirror image. We emphasize that the hierarchical sets under consideration here consist of self-similar Descartes configuration that share the boundary circle. Two mutually tangent circles  of curvatures $\kappa_1$ and $\kappa_2$  are respectively the boundary circles for the red and the green hierarchy and all the circles in the Pappus chains are tangent to both these boundary circles.

 \begin{figure}[htbp] 
  \includegraphics[width = .8 \linewidth,height=0.3 \linewidth]{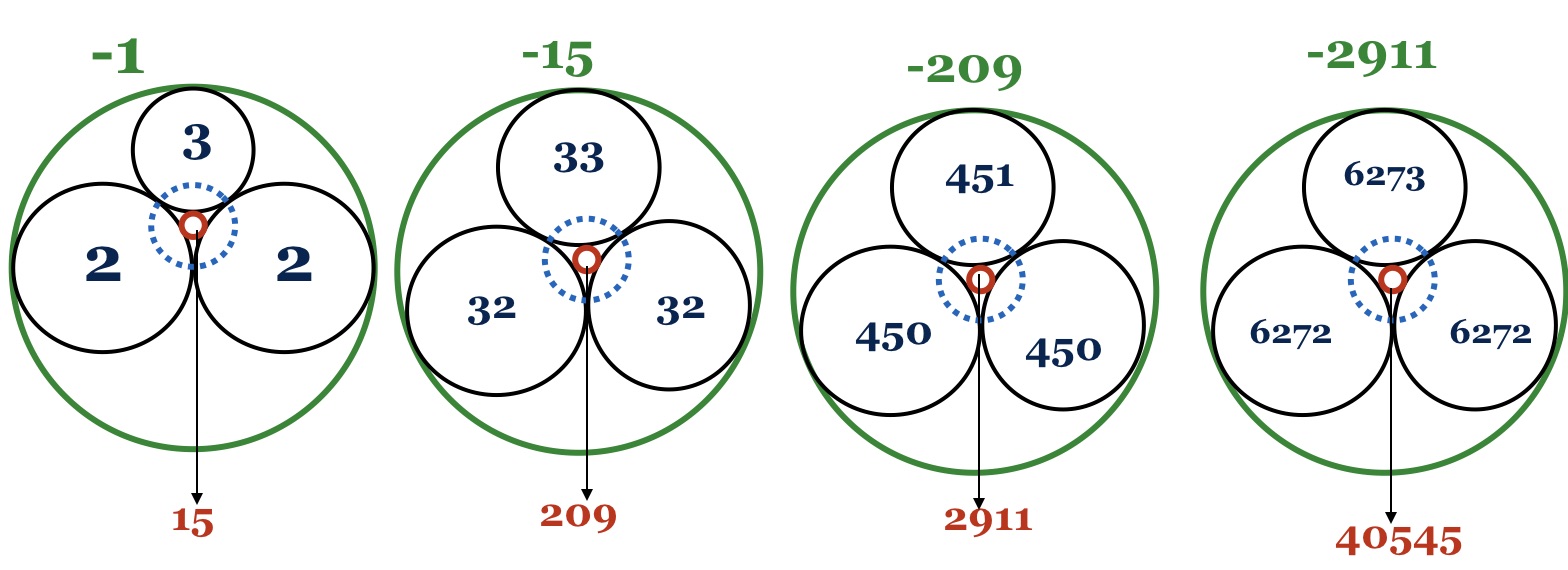}
\leavevmode \caption{ Self-similar kaleidoscopic symmetry in a hierarchy where the outer circle (green) becomes the inner circle (red) at the next level of the hierarchy. Inner and outer circles are mirror images through (dotted) blue circle. As we zoom into the innermost circles (  whose curvatures are shown
 below the arrows), we see a scaled version of the original pattern. The  sequence of ratios of curvatures  $( 15/1,209/15, 2911/209, 40545/2911, ...)$ approaches a constant $(2+\sqrt{3})^2$ . }
\label{SSM}
\end{figure}

 A rather exotic kind of hierarchical kaleidoscopic phenomena prevalent in packing of circles is shown in Figure (\ref{SSM}) where the object and its mirror image are self-similar. In other words, there exists 
some special set of recursive patterns in the gasket that are captured by  kaleidoscopes that reproduce the exact replica of the original set after appropriate scaling. Relationship of this self-similar kaleidoscope that exhibits three-fold symmetry with the gasket in figure ~ (\ref{iagall}) will be discussed later.

\section{ Ford Circles}
\begin{figure}[htbp] 
 \includegraphics[width = .6 \linewidth,height=0.35 \linewidth]{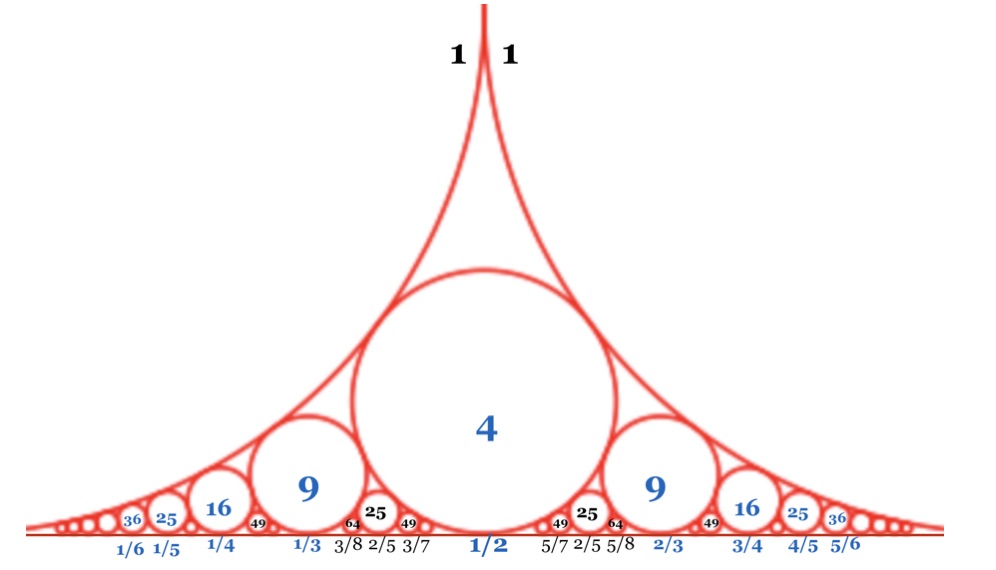} 
\leavevmode \caption{ Gasket made up of Ford circles where each circle with center at $(\frac{p}{q}, \frac{1}{2q^2})$ represents a primitive fraction $\frac{p}{q}$. After scaling by a factor of $2$, curvatures of all Ford circles are integer-squared. }
\label{ford}
\end{figure}
  Perhaps the simplest example of an hierarchical set consisting of packing of circles is a gasket made up circles  where the curvatures of all the circles are integer-squared. These circles are tangent to  the $x$-axis and  provide a pictorial representation of fractions as discovered by an American mathematician Lester Ford 
in  $1938$. Ford showed\cite{Ford} that at each rational point $\frac{p}{q}$ where $p$ and $q$ are relatively prime, one can 
draw a circle of radius $\frac{1}{2q^2}$ and whose center is the point $(x, y) = (\frac{p}{q}, \frac{1}{2q^2}) $. Tangent to the x-axis in the upper half of the $xy$-plane, the curvatures of the Ford circles are $2q^2$. Circles in the figure (\ref{ford}) have their curvatures $q^2$, scaled by a factor of half from Ford circles. We will  continue to refer them as
the Ford circle representing the fraction $\frac{p}{q}$ and  label them as $\kappa_{\frac{p}{q}}$.

  The key characteristic of the Ford circles is the fact that two Ford circles representing two distinct fractions {\it never} intersect. The closest they can come  is being tangent to each other. This happens when the two fractions are ``Farey neighbors", neighboring fractions in the Farey tree. Three mutually tangent Ford circles representing three  (left, center and right ) fractions $\phi_L = \frac{p_L}{q_L}$, $\phi_c = \frac{p_c}{q_c}$ , $\phi_R = \frac{p_R}{q_R}$  , all tangent to the $x$-axis form a special case of Descartes configuration obeying the Farey relation\cite{Ford},
  \begin{equation}
  \frac{p_c}{q_c} = \frac{p_L+ p_R}{q_L+q_R}
  \label{FR}
  \end{equation}
 Consequently, the triplet of Farey neighbors, referred as the {\it friendly triplet}   
  satisfy the follow identities:

\begin{equation}
|q_{\rm L} p_{\rm R} - q_{\rm R} p_{\rm L} | = 1,\,\,\ |q_{\rm L} p_{\rm c} - q_{\rm c} p_{\rm L} | 
= 1,\,\,\ |q_{\rm R} p_{\rm c} - q_{\rm c} p_{\rm R} | = 1
\ .
\label{FR}
\end{equation}

Figure ~(\ref{ford}) shows an Apollonian gasket made up of the Ford circles where the curvatures are integer squared and we will refer such gaskets as the `` Ford Apollonian gasket". The dual of the Ford Apollonian gasket a subset of circles
seen in  figure (\ref{iagall}) that exhibit reflection symmetry about the $x$-axis satisfying $\kappa_1+\kappa_2= \kappa_3+\kappa_4$ .
In the rest of the paper, we will denote Descartes configuration of Ford circles  representing the friendly  triplet $[\frac{p_L}{q_L},\frac{p_c}{q_c},\frac{p_R}{q_R}]$ as $( \kappa_c, \kappa_R, \kappa_L)$  and its dual as $( -\kappa_0, \kappa_1, \kappa_2, \kappa_3)$. Here $ \kappa_c > \kappa_R > \kappa_L$ and
$ \kappa_0 < \kappa_1 < \kappa_2 < \kappa_3$. 

It was pointed out by Richard Friedberg\cite{RF} that the Descartes’s theorem when applied to Ford circles is a consequence of a property of homogeneous functions of
any three variables $(\lambda_1, \lambda_2, \lambda_3)$ with  $\lambda_1 \pm \lambda_2 \pm \lambda_3=0$. In view of the Farey relation $ q_c = q_L + q_R$, as shown in the  Appendix B, this gives $2(q_1^4+q_2^4+q_3^2) = ( q^2_1+q_2^2+q_3^2)^2$ which is the simplified version ( $ \kappa_4= q_4^2 =0$ ) of the Descartes theorem with three mutually tangent Ford circles, tangent to the $x$-axis.

The Ford Apollonian gasket  has certain properties that  play an important role in determining the self-similar scalings of the entire Apollonian gasket even through the Ford circles form only a small subset of the  fractal. Central result of this paper is that 
packing of the Ford circles {\it determines} the self-similar properties of the entire gasket. 
Further, there is a dichotomy in the recursive structure of the gasket as half of the gasket seems to follow the recursive pattern of the Pythagorean tree\cite{SatPT}. In addition,
 Ford circles also have a partner - a symmetric Descartes configuration where two of the three inner circles have same curvature and there exists a very special configuration that evolves into three-fold symmetry\cite{book}  as shown in Fig. (\ref{SSM}). These three-fold symmetric configurations that exhibit a self-similar hierarchy with Kaleidoscopic symmetries lay scattered throughout the whole gasket.

\section{ Self-Similarity of Ford Apollonian Gasket}
   \begin{figure}[htbp] 
\includegraphics[width = .6 \linewidth,height=.4 \linewidth]{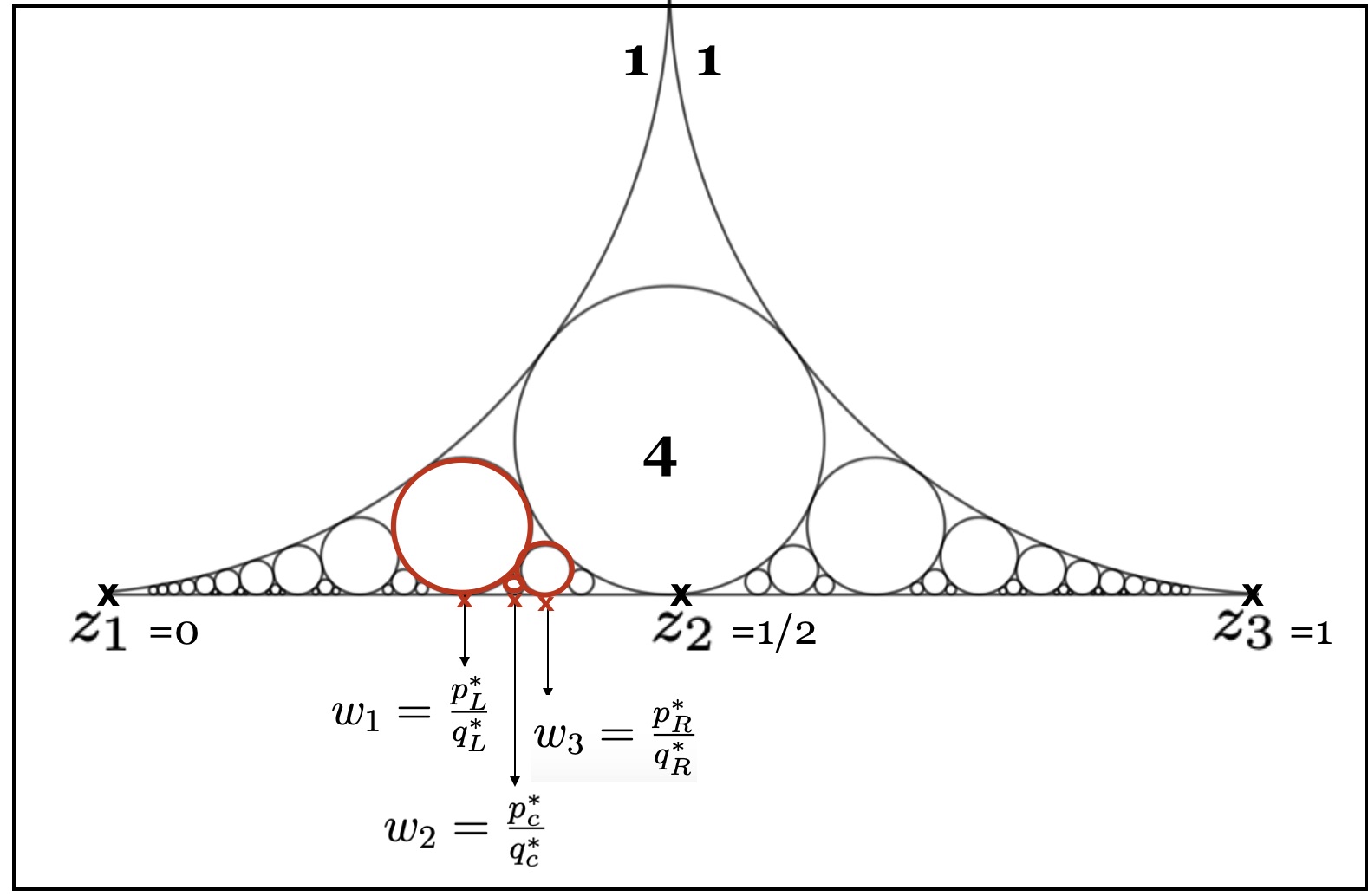}
\leavevmode \caption{With the  root configuration $(\kappa_c, \kappa_R, \kappa_L)=(4,1,1)$ representing three fractions $(z_1,z_2,z_3)=(\frac{0}{1}, \frac{1}{2},\frac{1}{1})$, we pick any other configuration of three mutually tangent Ford circles labeling three fractions as
$( w_1, w_2, w_3)=(\frac{p^*_L}{q^*_L}, \frac{p^*_c}{q^*_c},\frac{p^*_R}{q^*_R})$ representing left, central and right circle.  In principle, one can choose any of the three tangency points  and their corresponding image of the Descartes configurations, that is, $z_i$ and $w_i$ need not lie on the $x$-axis. The choice
indicated here simplifies the derivation of the recursion relation.}
\label{FC}
\end{figure}

 To establish recursion relations associated with the nesting of the Ford circles, we relate  two  sets of  triplets: namely
 $(z_1, z_2, z_3)$ and $(w_1, w_2, w_3)$ as shown in figure ~(\ref{FC}). The self-similar hierarchy preserves this relationship at all levels and we will denote two successive  sets as $[\phi_L(l), \phi_c(l), \phi_R(l)]$
and $[\phi_L(l+1), \phi_c(l+1), \phi_R(l+1)]$.
The recursion relation underlying the self-similar pattern is a fixed M\"{o}bius transformation  determined by mapping the triplets
  $( z_1, z_2, z_3)=( 0, \frac{1}{2}, 1)$ to the triple
 $( w_1,w_2, w_3) = (\frac{p^*_L}{q^*_L}, \frac{p^*_c}{q^*_c},\frac{p^*_R}{q^*_R}) =( \phi^*_L(1), \phi^*_C(1), \phi^*_R(1) )$. 
 
 To determine the constants $(a,b,c,d)$  of the map $w = f(z) = \frac{az+b}{cz+d}$, we note that :
 \begin{equation}
 f(0)= \frac{b}{d}=\frac{p^*_L}{q^*_L},\,\,\,\,\ f(1/2)= \frac{a+b}{c+d}=\frac{p^*_R}{q^*_R},\,\,\,\,\ f(1)=\frac{a+b}{c+d}=\frac{p^*_R}{q^*_R}
 \end{equation}
 Using the  Farey relation (Eq. (\ref{FR})), we get
 \begin{equation}
 a=p^*_R-p^*_L,  \,\ b=p^*_L, \,\ c = q^*_R-q^*_L  ,  \, \,\ d = q^*_L
 \label{abcd1}
 \end{equation} 
 
 Therefore, the  M\"{o}bius map that underlies the recursive structure of the Ford Apollonian gasket is given by a matrix, which we denote as $\mathcal{F^*}$.
The corresponding $\phi= \frac{p}{q}$ recursions  can be written as,
\begin{equation}
 \phi(l+1) =  \frac{ (p^*_R-p^*_L) \phi (l) + p^*_L}{ (q^*_R-q^*_L) \phi (l) + q^*_L} \equiv \frac{ap+bq}{cp+dq} \equiv  \mathcal{F}^*
 \label{Mmap}
 \end{equation}
 This  can also be written in terms of the recursions for  numerator $p_x$ and the denominator $q_x$ of the fraction $\frac{p_x}{q_x}$  where $x=L, c, R$, 
\begin{equation}
\left(\begin{array}{c}
p_x(l+1)\\
q_x(l+1)
\end{array}\right)   =   \mathcal{F}^*  \left(\begin{array}{c}
p_x(l) \\
q_x(l)
\end{array}\right) 
\end{equation}

    \begin{figure}[htbp] 
\includegraphics[width = .7 \linewidth,height=.4 \linewidth]{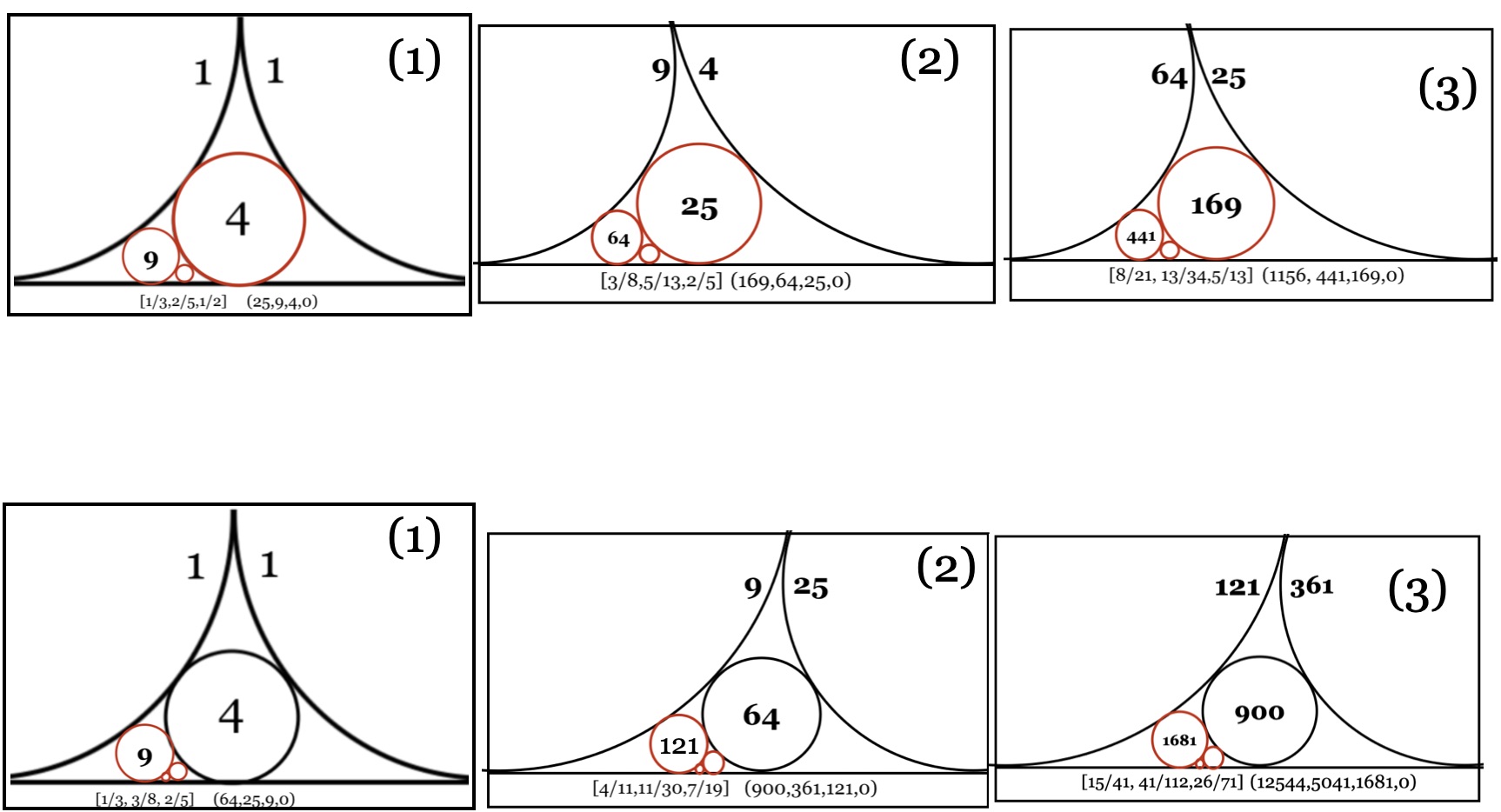}
\leavevmode \caption{ Upper and lower panels respectively illustrate $n^*=1$ and $n^*=2$ self-similar hierarchies where three panels in each case represent three levels of the hierarchy consisting of Ford circles.  Starting from the root  Apollonian( black circles), we pick an Apollonian ( red circles ).The self-similar hierarchy corresponds to continuing this process iteratively where the relationship between two successive iterations is preserved. That is the relation between $(4,1,1)$ and $( 25,9,4)$ is the same as the relation between $( 25, 9, 4)$ and $( 169, 64, 25)$ and so on. Asymptotic scale invariance
 is signaled as $25/4=6.25, 169/25=6.76 \rightarrow (\frac{3+\sqrt{5}}{2})^2 \approx 6.85$ ( golden-hierarchy, upper panels) and
 $ 121/9=13.44, 1681/121=13.89 \rightarrow ( 2+\sqrt{3})^2 \approx 13.93$ ( diamond hierarchy, lower panels).  As described later, the lower hierarchy  where the parity of $q_c$ is conserved, mimics the tree structure of the Pythagorean tree.}
\label{SS}
\end{figure}
 
\section{Scaling}

With unit determinant, the eigenvalues of the transformation matrix  $\mathcal{F^*}$ are real. Denoting the pair of these eigenvalues as $(\zeta, \zeta^{-1})$, $\zeta > 1$ determines\cite{book}  the asymptotic scaling factors as:

\begin{equation}
\zeta= \lim _{ l \rightarrow \infty}  \frac{p_x(l+1)}{p_x(l)} = \lim _{ l \rightarrow \infty}  \frac{q_x(l+1)}{q_x(l)}  =\frac{( q^*_L+ p^*_R-p^*_L)}{2}\pm\sqrt{\left(\frac{q^*_ L+  p^*_R-p^*_L}{2}\right)^2-1} ,
 \label{zeta1}
 \end{equation}
 The ratio of curvatures of the circle at two successive levels is given by,
 \begin{equation}
 \lim_{ l \rightarrow \infty}  \frac{\kappa (l+1)}{ \kappa(l)}  \rightarrow \zeta^2
  \end{equation}
Expressed as a  continued fraction expansion, these irrationals  satisfy the quadratic equation $\zeta^2+n^*\zeta-n^*$ whose solutions are given by,
 \begin{equation}
 \zeta= [n^*+1; \overline{ 1, n^*}], \,\ n^* = q^*_L + p^*_R-p^*_L-2
 \label{star}
 \end{equation}
  Figure (\ref{SS}) shows two examples of self-similar Descartes configurations of Ford circles corresponding to $n^*=1$ (upper panel)  and $n^*=2$ (lower panel). Respectively referred as the {\it golden} and  the {\it diamond} hierarchies\cite{book}, they correspond to the  M\"{o}bius transformations $f(z)=\frac{1}{-z+3} \equiv \left( \begin{array}{cc} 0  & 1  \\  -1 & 3   \\  \end{array}\right)$ and $f(z)=\frac{z+1}{2z+3}  \equiv \left( \begin{array}{cc} 1  & 1  \\  2 & 3   \\  \end{array}\right) $. The eigenvalues of these two matrices  $\frac{3\pm \sqrt{5}}{2}$ and $2\pm \sqrt{3}$ determine the scaling factor $\zeta$. We emphasize that the M\"{o}bius map (\ref{Mmap}) captures the self-similarity of the whole gasket provided the root configuration is $(4,1,1)$. General case with arbitrary root is described later.

   \section{ Self-Similar Hierarchies with non-Ford circles }
 
\begin{figure}[htbp] 
\includegraphics[width = .68 \linewidth,height=.33 \linewidth]{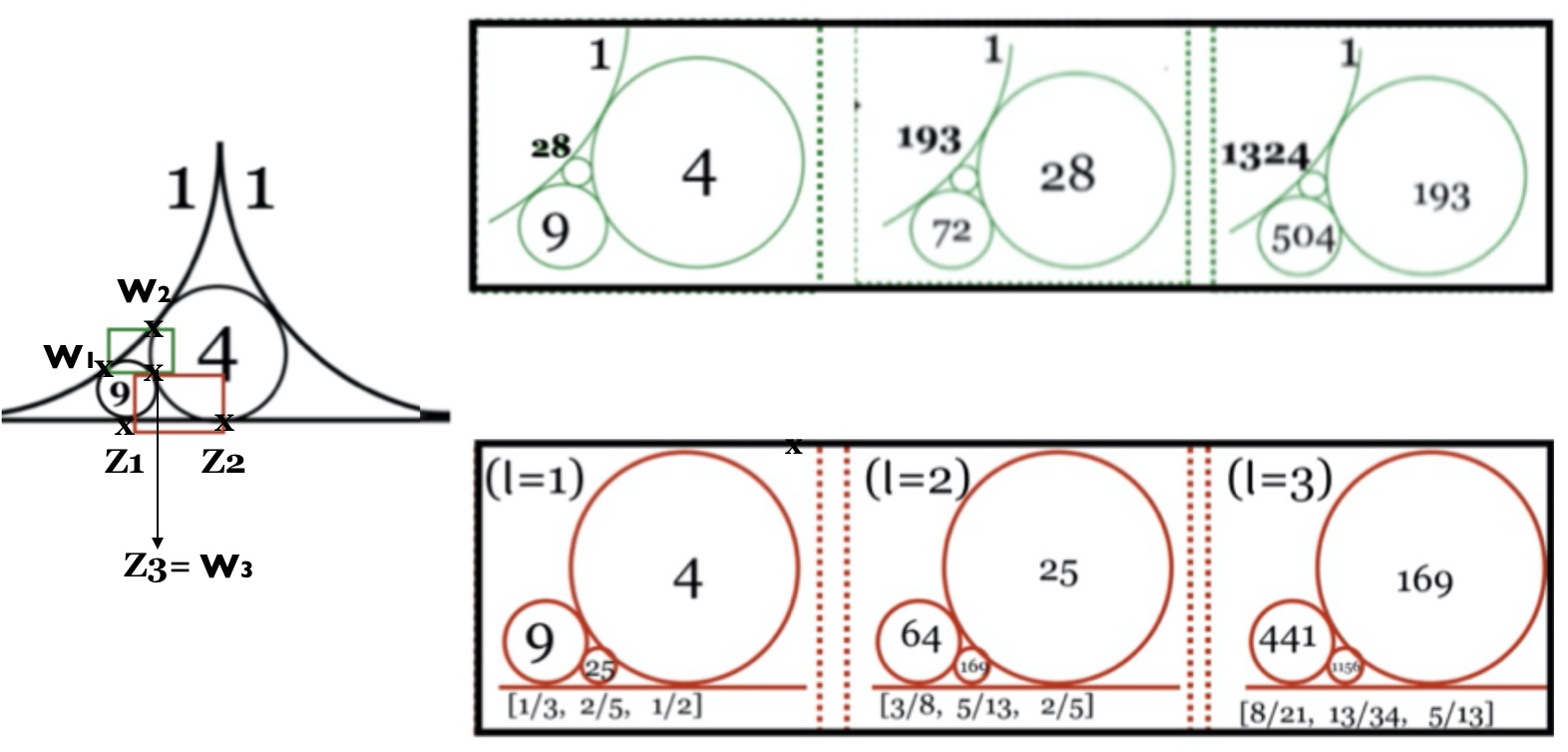}
\leavevmode \caption{ The green and the red boxes in the left panel respectively show the locations of the nesting patterns  of the green and the red circles -  the two hierarchies  shown in the upper and lower right panels. These two hierarchies  exhibit same scaling exponents. There distinct Descartes configurations of red and green circles show three levels of  the self-similar hierarchies. }
\label{E2C}
\end{figure}

 The M\"{o}bius transformations map circles to circles, preserving tangencies. Therefore,  every Descartes configurations consisting of  ``non-Ford circles" -  circles that are not tangent to $x$-axis,  can be related to the Ford circles. The process of describing the self-similarity of the non-Ford hierarchies  involves  first finding the conformal image of the hierarchy of the non-Ford circles 
 to  the corresponding Ford circles. One can then use the equation (\ref{Mmap}) to describe the self-similar characteristics.  Fig. (\ref{E2C}) illustrates this process when the non-Ford circle hierarchical pattern is tangent to $\kappa_{\frac{0}{1}}$. In this example, we first obtain
the M\"{o}bius transformation that maps  the $x$-axis to $\kappa_{\frac{0}{1}}$  in the complex plane.  As shown in the left panel,
the three points on the $x$-axis and their conformal images on $\kappa_{\frac{0}{1}}$ can be chosen as:
$ z_1=0 \rightarrow w_1 = 0,\,\ ,\ z_2 = \frac{1}{2} \rightarrow  w_2=\frac{2}{5}+ \frac{i}{5}, \,\,\,\  z_3 = \frac{1}{3} \rightarrow w_3 = \frac{3}{10}+\frac{i}{10}$
  The resulting map   is given by
$\mathcal{B}: z \rightarrow w=f(z)= \frac{z}{-i z+1 }$. Its inverse transforms the green circles to the red circles, and therefore, 
 the hierarchical structure of the non-Ford circles is given by the general map 
$\mathcal{B}^{-1} \mathcal{F^*} \mathcal{B}$.  Eq.  (\ref{Mmap}) determines $F^*$ where $p^*_L= 1$, $q^*_L = 3$, $p_R^*=1$ and $q_R^*=2$.
 
\begin{figure}[htbp] 
\includegraphics[width = .65 \linewidth,height=.3 \linewidth]{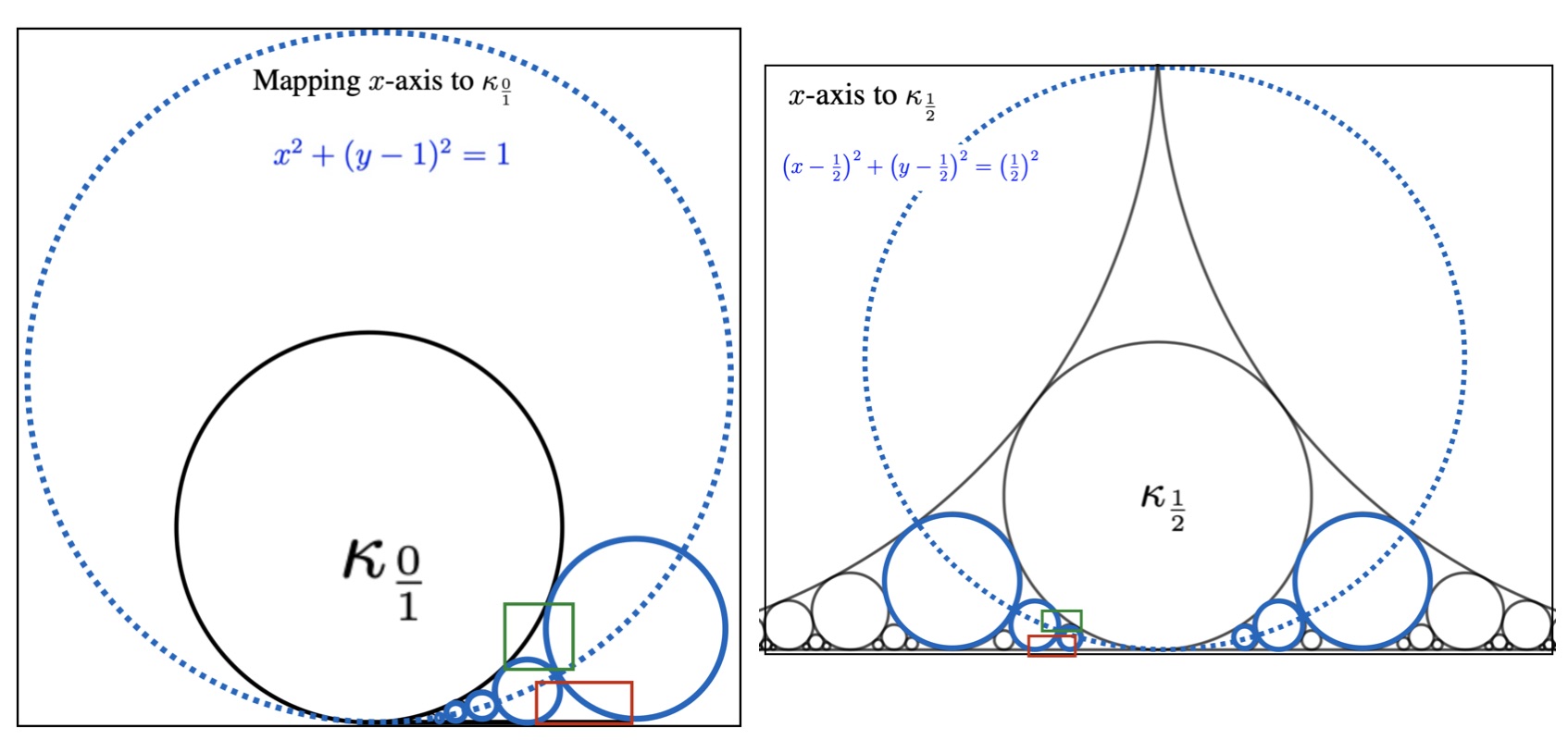}
\leavevmode \caption{ Left and right panels  show the Pappus chains ( blue circles )  and the corresponding Pappus mirrors ( dotted blue circles )  that reflect the x-axis to $\kappa_{\frac{0}{1}}$ and $\kappa_{\frac{1}{2}}$ respectively. The red and the green boxes are the triangular regions hosting circles that are mirror images of each other. }
\label{K2}
\end{figure}

 Table (1) gives some examples of  M\"{o}bius transformations that relate $x$-axis to other circles in the $\mathcal{IAG}$.  The corresponding ``mirrors" listed in the table ( see Fig. (\ref{K2}) ) provides geometrical 
 representation of this mapping. We note that  with the first three entries in the Table that show the mappings of the $x$-axis to  
the Ford circles $\kappa_{\frac{p}{q}}$ 
 can be written as  $g(z)= \frac{1}{q^2} f( q^2(z-\frac{p}{q})) + \frac{p}{q}$ where $f(z) = \frac{z}{-iz+1}$ maps the $x$-axis to $\kappa_{\frac{0}{1}}$. These transformations form a subgroup of $SL(2,C)$ with trace,  $a+d=2$. They  are ``parabolic as they have only one fixed point. 
 
All the transformations described in the table
have unit determinant and real trace.  Such transformations, which we denote as
 $\mathcal{B}$,  that relate boundary circle $\kappa_{\frac{p}{q}}$ of an hierarchy to the $x$-axis, the self-similar recursive pattern is given by $\mathcal{B} \mathcal{F}^* \mathcal{B}^{-1}$. Therefore, the scaling factor associated with the hierarchy is  same as the scaling factor associated with their corresponding Ford circle hierarchy.
 
 \begin{table}
\begin{tabular}{| c | c |  c  |}
\hline
Nature of Map & 
f(z) \,\, &  Equation of the Mirror   \\ 
\hline
$x$- axis $\rightarrow \kappa_{\frac{0}{1}}$ \,\,  & $ \frac{z}{-iz+1}$ \,\, & $x^2+(y-1)^2=1$  \\ 
\hline
$x$- axis $\rightarrow \kappa_{\frac{1}{2}}$ \,\,  & $  \frac{(1-2i)z+i}{-4iz+(1+2i)}$ \,\, & $(x-\frac{1}{2})^2+(y-\frac{1}{2})^2=(\frac{1}{2})^2$  \\ 
\hline
$x$- axis $\rightarrow \kappa_{\frac{1}{3}}$ \,\,  & $ \frac{(1-3i)z+i}{-9iz+(1+3i)}$ \,\, & $(x-\frac{1}{3})^2+(y-\frac{2}{9})^2=(\frac{2}{9})^2 $  \\ 
\hline
$x$- axis $\rightarrow \kappa_{\bar{\frac{2}{5}}}$ \,\,  & $\frac{(2-5i)z+2i}{-14i z+(2+5i)}$ \,\, & $(x-\frac{5}{14})^2+(y-\frac{2}{7})^2=(\frac{1}{7})^2 $   \\ 
\hline
\end{tabular}
\caption{ Examples of M\"{o}bius transformations and the corresponding mirrors that map $x$-axis to the Ford circles $\kappa_{\frac{p}{q}}$ or the conformal images of the Ford circle, tangent to $\kappa_{\frac{0}{1}}$, which we denote as  $\kappa_{\bar{\frac{p}{q}}}$.  Fig. (\ref{K2}) shows two examples of the mirrors, corresponding to the first two entries in this table.}.
\label{T1}
\end{table}

\begin{figure}[htbp] 
\includegraphics[width = .38 \linewidth,height=.5 \linewidth]{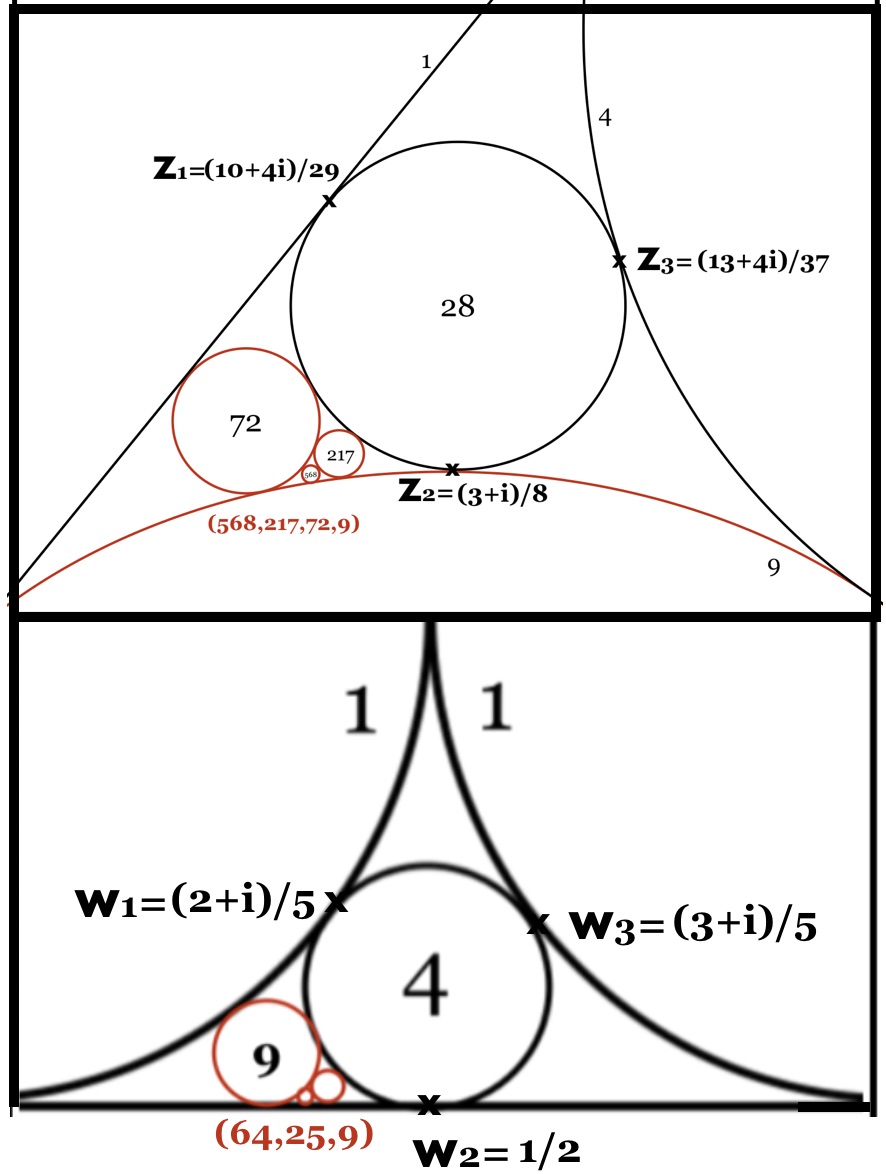}
\leavevmode \caption{ (Upper) Starting with a root Apollonian $(28,9,4,1)$, we select a configuration $(568,217,72,9)$ (red circles) that bears the same relationship with the root as  $(64,25,9)$ (lower panel) bears with the main root $(4,1,1)$ associated with the Ford Apollonian gasket. Using $(z_1, z_2, z_3)$ from the upper graph and their corresponding $(w_1, w_2, w_3)$ in the lower graph, 
we obtain the M\"{o}bius map $\mathcal{B}: w =\frac{(-9+6i)z+(4-i)}{(-14+16i)z+(7-4i)}$  that transforms the root $(28,9,4,1)$ to $(4,1,1)$.}
\label{C1E1}
\end{figure}

Finally, when the root configuration is not  the  $(4,1,1)$, the chosen root is conformally mapped to $(4,1,1)$ as illustrated in figure ~(\ref{C1E1}).
In this example, using the triplet $(z_1, z_2, z_3)$ and its image $(w_1, w_2, w_3)$ 
leads to the map $\mathcal{B}:  w =f(z)=\frac{(3-i)z-1}{(-3i)z+i}$ that transforms the chosen root to $(4,1,1)$. We then use the map $ \mathcal{F^*}: z \rightarrow \frac{z+1}{2z+3}$  that describes the self-similar recursions in the Ford circles in the lower panel. In other words, the self-similar characteristics of the upper panel hierarchy  is described by the map $\mathcal{B} \mathcal{F}^* \mathcal{B}^{-1}$.

\section{ Self-Similar Kaleidoscope}

\begin{figure}[htbp] 
\includegraphics[width = .63 \linewidth,height=0.3 \linewidth]{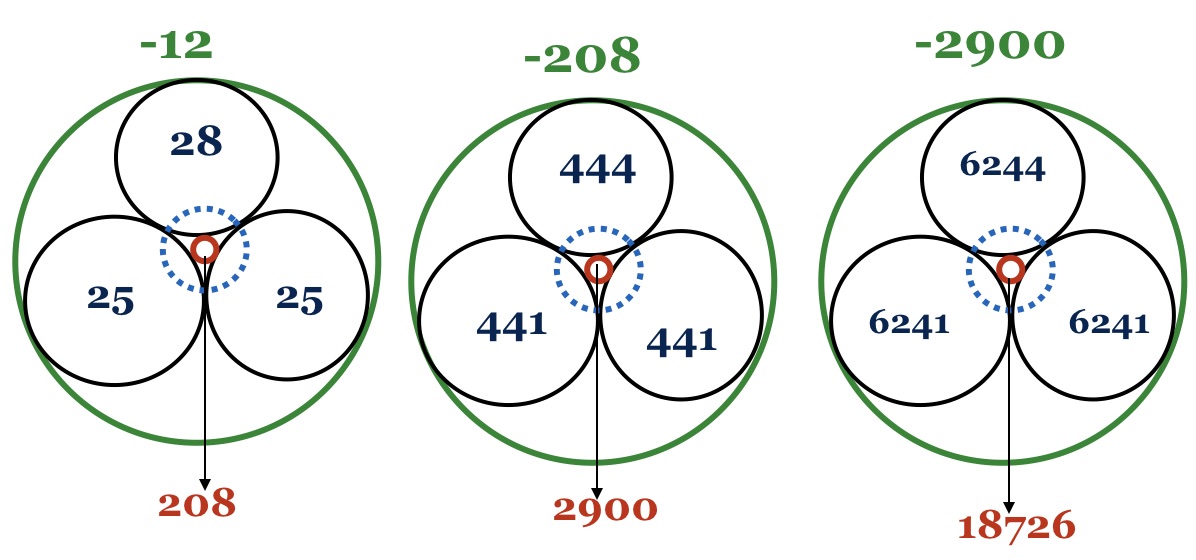} 
\leavevmode \caption{ Symmetric partners corresponding to lower panel of Fig. (\ref{E2C}) where $\Delta=3$. The corresponding Fig. (\ref{SSM}) describes the symmetric partners of the upper panel of the Fig. (\ref{E2C})
corresponding to $\Delta=1$. Differing in the initial root, 
these two hierarchies belong to same universality class characterized by scaling exponent $\zeta=( 2+ \sqrt{3})$.}
\label{fordR}
\end{figure}

Among the infinity of self-similar hierarchies in the $\mathcal{IAG}$, there is one type of a hierarchy that is characterized by a very unique symmetry shown in Fig. (\ref{SSM}) and also in Fig. (\ref{fordR}).
These configurations  are
related to the Ford Apollonian as
every Descartes configuration of the Ford circles can  be mapped to another Descartes configuration of the non-Ford circles where two of the inner circles have the same curvatures\cite{book}.  In other words, 
given $( \kappa_c, \kappa_R, \kappa_L)$  or its dual $(- \kappa_0, \kappa_1, \kappa_2, \kappa_3)$, we can find a Descartes configuration which we denote as
 $( \kappa^s_a, \kappa^s_b, \kappa^s_b,  \kappa^s_c)$ where,
 
\begin{equation}
\kappa_a^s =  \eta  \kappa_0, \,\,\  \kappa_b^s =  \frac{\eta}{2} ( \kappa_1+ \kappa_2) ,\,\,\ \kappa_c^s = \eta (2\kappa_1 - \kappa_0)
\end{equation}
\begin{equation}
\kappa_a^s = \frac{\eta}{2} ( \kappa_c - \kappa_R -\kappa_L),\,\,\ \kappa_b^s =  \frac{\eta}{2} \kappa_c, \,\,\, \kappa_c^s = \frac{\eta}{2} ( \kappa_c - \kappa_R + 3 \kappa_L) 
\end{equation}

Here $\eta=1$ for $\kappa_0$-odd and $\eta =2$ for $\kappa_0$-even.

Among these symmetric configurations, there exists a class of hierarchies where the above configuration evolves towards an asymptotic three-fold symmetry\cite{book} as shown in figures (\ref{SSM}) and (\ref{fordR}).This is due to an invariant $\Delta$\cite{book},
\begin{equation}
 \Delta=\kappa^s_b-\kappa^s_c, 
 \label{del}
 \end{equation}
which remains unchanged under iteration of the recursions. In our detailed study of various such configurations, differing in $\Delta$, which was always found to be a prime number, all such hierarchies
 belong to  the nested configurations of the diamond hierarchy.  Relationship of three fold symmetry with diamond hierarchy follows from Eq. (\ref{Q4})\cite{book} as with  $\kappa_1=\kappa_2=\kappa_3$,  $\frac{\kappa_4}{\hat{\kappa}_4} = ( 2+\sqrt{3})^2$. 

\section{ Apollonian-Pythagorean Meet }

 \begin{figure}[htbp] 
\includegraphics[width = .75\linewidth,height=0.7\linewidth]{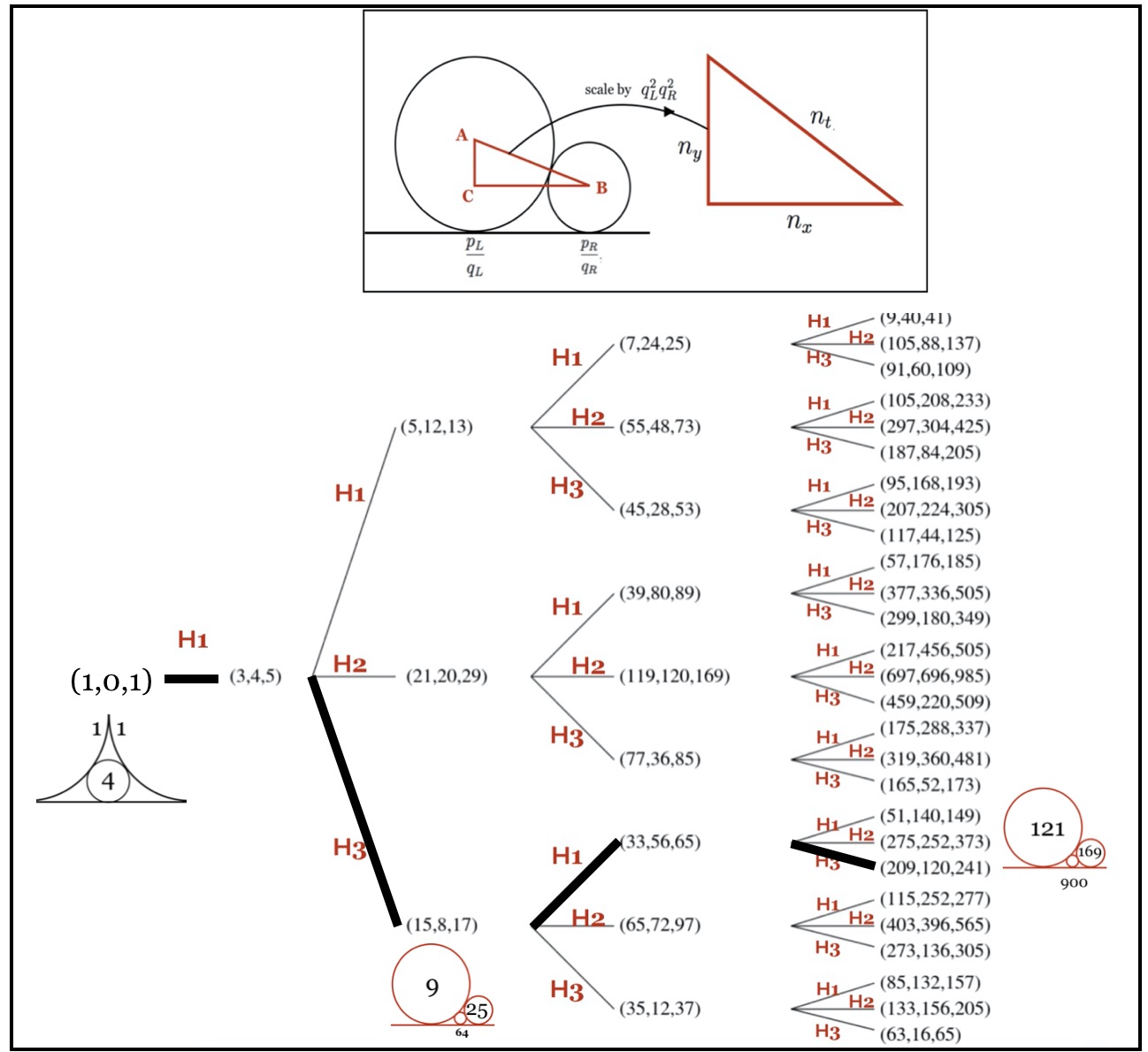} 
\leavevmode \caption{ Upper panel shows mapping between two mutually tangent Ford circles to a Pythagorean triplet\cite{SatPT}. The lower panel shows the diamond hierarchy described by $H_1H_3$, highlighted  with thick dark lines along with the Descartes configurations. Starting with the root $(4,1,1)$ that corresponds to the Pythagorean triplet $(1,0,1)$, the triplets $(15,8,17)$ and $(209, 120, 241)$ represent the level 1 and level 2 of Descartes configurations shown in red circles.}
\label{PT}
\end{figure}

The self-similar recursions described by the Ford Apollonian gasket $( \kappa_c, \kappa_R, \kappa_L)$ are found to belong to two classes depending upon whether they preserve ``parity"-  that is even or oddness of $\kappa_c$ (or $q_c$). 
It turns out that this parity conserving feature is directly correlated with the parity of 
the  integer $n^*$ ( see Eq. (\ref{star}) ). For  the recursions characterized by even $n^*$, the parity is conserved while for  the recursions with odd-$n^*$,  the parity is not conserved. Fig. (\ref{SS}) shows examples of both of these cases.
 
 An interesting aspect of the parity conserving recursions is that they are described by the Pythagorean treecite{Hall, J3, SatPT} . As shown in Fig. (\ref{PT}),  the Pythagorean tree  provides orderly
arrangements of all primitive Pythagorean triplets $(n_x, n_y, n_t)$ where $ n_x^2 + n_y^2 = n_t^2$\ using three  matrices $(H_1, H_2, H_3)$:

\begin{equation}
H_1  =   \left( \begin{array}{ccc} 1 & -2 & 2  \\    2 & -1 & 2  \\  2 & -2  & 3 \\ \end{array}\right),\,\,\
H_2  =  \left( \begin{array}{ccc} 1 & 2 & 2  \\    2 & 1 & 2  \\  2 & 2  & 3 \\ \end{array}\right) ,\,\,\,\ 
H_3  = \left( \begin{array}{ccc} -1 & 2 & 2  \\    -2 & 1 & 2  \\  -2 & 2  & 3 \\ \end{array}\right) 
\label{H3}
\end{equation}

For example, with Pythagorean triple $ v = ( 1, 0,1)$, we see that $H_1 v ^T=(5,12,13)^T$, and so on.
There is also a two-dimensional representation of $H_i$
using $2000$ year old Euclid parametrization of the Pythagorean triplets  given in terms of pair of integers $(q_R,q_L)$, $ q_R > q_L$ where, 
\begin{equation}
 n_x =   \eta (q_L q_R), \,\  n_y = \eta ( \frac{q_R^2-q_L^2}{2}), \,\ n_z =\eta ( \frac{q_R^2+q_L^2}{2}). 
 \end{equation}
Here $\eta=1$  when $q_c$ is even and hence  $q_L$ and $q_R$ are both odd integers and $\eta=2$ otherwise. 
The possibilities of generating  three  more triplets from a given triplet, expressed in terms a  pair $(q_R,q_L)$ can be written as three 
 $ 2\times 2$ matrices,

\begin{equation}
h_1 =  \left( \begin{array}{cc} 1   & 2  \\  0 & 1   \\  \end{array}\right), \quad  h_2 =  \left( \begin{array}{cc} 2   & 1  \\  1 & 0   \\  \end{array}\right),\,\,\,\ h_3 =  \left( \begin{array}{cc} 2   & -1  \\  1 & 0   \\  \end{array}\right)
\label{h2}
\end{equation}

 The hierarchical pattern of  parity-conserving Ford Apollonian can be associated  with the hierarchical pattern of the Pythagorean tree\cite{SatPT}. 
Given a Ford Apollonian  $(\kappa_c, \kappa_R, \kappa_L)$, a linear transformation relates the Apollonian to a Pythagorean triplet $(n_x, n_y, n_t)$ as:

\begin{eqnarray}
\left( \begin{array}{c} n_x\\ n_y \\ n_t  \\ \end{array}\right) & = & \frac{1}{\eta} \left( \begin{array}{ccc} 1 & -1 & -1  \\    0 & 1 & -1  \\  0 & 1  & 1 \\ \end{array}\right)  
  \left( \begin{array}{c}  \kappa_c \\ \kappa_R\\ \kappa_L  \\ \end{array}\right) ,
  \label{b2pt}
\end{eqnarray}

With $\eta=1$ or $2$,  there are two Pythagorean trees related simply by swapping $n_x$ and $n_y$.The key point to be noted here is that the Pythagorean triplets in each of these two trees retain their parity along the tree. Consequently, only parity conserving recursions of the Ford Apollonian are described by the Pythagorean tree.
As a consequence of this dichotomy in characterization of self-similarity of an Apollonian gasket, half of the  hierarchies in the  Ford Apollonian gaskets share the recursive structure of the Pythagorean tree. The even parity
hierarchies are generated by $(3,4,5)$ using $(H_1, H_2, H_3)$ while the odd-parity hierarchy is generated by $(4,3,5)$ using $( H_3, H_2, H_1)$. The  diamond hierarchies which conserves $q_c$ parity as described above correspond to $H_3H_1$. The hierarchical aspect of the golden hierarchy is not described by the tree as the corresponding recursions do not conserve parity as seen in figure (\ref{SS}). Recursions correspond to zigzagging between the $\eta=1$ and $\eta=2$
Pythagorean trees.
  
As a generalization of the relation between the Pythagorean triplet and the Descartes configuration of the Ford circles, we note that
the four curvatures  $(\kappa_1, \kappa_2, \kappa_3, \kappa_4)$ of any Descartes configuration are related to the well known quadruplets - $( N_x, N_y, N_z, N_t)$\cite{J4} , known as
the  Lorentz quadruplets,  where
 $ N_x^2+N_y^2+N_z^2=N_t^2$\cite{J2}. This relationship is given by,   

\begin{eqnarray}
 \left( \begin{array}{c} N_x  \\    N_y  \\  N_z  \\  N_t   \\ \end{array}\right) =  \left( \begin{array}{cccc} 1 & -1 & -1  & -1 \\    0 & 0 & 0 & 2  \\  0 & 1 & -1 & 0\\ 1 & 1 & 2& 1  \\ \end{array}\right) 
  \left( \begin{array}{c} \kappa_1  \\  \kappa_2  \\  \kappa_3 \\  \kappa_4 \\ \end{array}\right) . \end{eqnarray} 
 
 It turns out that the Lorentz quadruplets associated with the four curvatures  need not be primitive as is found to be case with the golden hierarchy. Furthermore, to best of our knowledge, there is no known quadruplet tree that generates all the primitive or non-primitive  Lorentz quadruplets. Therefore, the representation of the $\mathcal{IAG}$ with a  tree like structure of Lorentz quadruplets remains an open question.
 
 \section{ From simple Geometry to Elegant Group Theory}
 
The geometric construction of kaleidoscopic images to describe Apollonian packing can be formulated using Apollonian group where adding additional circles is accomplished by applying four matrices  to 
a root configuration\cite{AP}.  Eq. (\ref{l4})   determines the four  generators $S_i,  ( i = 1-4)$ of the group as,

 \begin{eqnarray*}
S_1 =  \left( \begin{array}{cccc} -1 & 2 & 2  & 2 \\    0 & 1 & 0 & 0  \\  0 & 0 & 1 & 0\\ 0 & 0 & 0& 1  \\ \end{array}\right),\
S_2= \left( \begin{array}{cccc} 1 & 0 & 0  & 0 \\  2 & -1 & 2 & 2  \\  0 & 0 & 1 & 0\\ 0 & 0 & 0& 1  \\ \end{array}\right),\
S_3  =  \left( \begin{array}{cccc} 1 & 0 & 0  & 0 \\    0 & 1 & 0 & 0  \\  2 & 2 & -1 & 2\\ 0 & 0 & 0& 1  \\ \end{array}\right),\
S_4= \left( \begin{array}{cccc} 1 & 0 & 0  & 0 \\    0 & 1 & 0 & 0  \\  0 & 0 & 1 & 0\\ 2 & 2 & 2& -1  \\ \end{array}\right)
 \label{S4}
\end{eqnarray*}

These generators do not describe self-similar characteristics of the gasket as the application of $S_i$ on a column matrix consisting of four curvatures in monotonic order does not result in a new set of curvatures  that preserves the monotonicity of the curvatures. That is, a self-similar hierarchy cannot be associated with a string of $S_i$.
However,  It is possible to define an alternative set of four matrices $D_i$ that preserves the monotonicity of the curvatures:
 
 \begin{eqnarray*}
D_1 =  \left( \begin{array}{cccc} 0 & 1 & 0  & 0 \\    0 & 0 & 1 & 0  \\  0 & 0 & 0 & 1\\ -1 & 2 & 2 & 2  \\ \end{array}\right) ,\
D_2= \left( \begin{array}{cccc} 2 & -1 & 2  & 2 \\  1 & 0 & 0 & 0  \\  0 & 0 & 1 & 0\\ 0 & 0 & 0& 1  \\ \end{array}\right),\
D_3  =  \left( \begin{array}{cccc} 2 & 2 & -1  & 2 \\    1 & 0 & 0 & 0  \\  0 & 1 & 0 & 0\\ 0 & 0 & 0& 1  \\ \end{array}\right),\ 
D_4= \left( \begin{array}{cccc} 2 & 2 & 2  & -1 \\    1 & 0 & 0 & 0  \\  0 & 1 & 0 & 0\\ 0 & 0 & 1& 0  \\ \end{array}\right)
 \label{D4}
\end{eqnarray*}

 \begin{figure}[htbp] 
\includegraphics[width = .7\linewidth,height=0.25\linewidth]{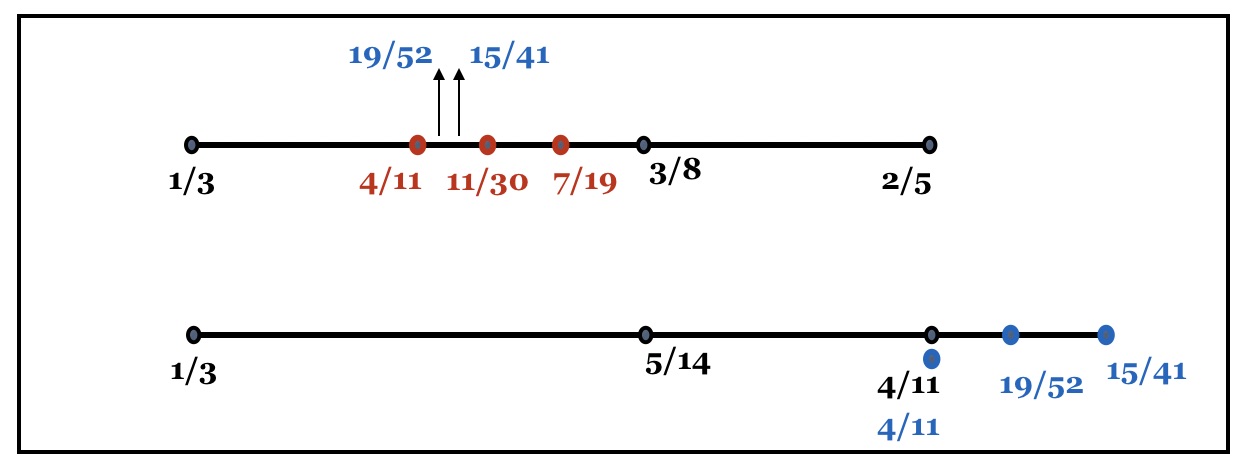} 
\leavevmode \caption{ Without showing the Ford circles, the sequences of upper (nested triplets of black and red dots) and lower (chain of black and blue dotes) fractions show two levels of two distinct hierarchies, corresponding to  Descartes configurations described by the friendly triplets $[0,1/2,1] \rightarrow [1/3, 3/8, 2/5] \rightarrow [4/11, 11/50, 7/19]...$ and
 $[ 0, 1/3, 1/4] \rightarrow [1/3, 5/14, 4/11] \rightarrow [ 4/11, 19/52, 15/51].......$. Recursions for both hierarchies are characterized by the same M\"{o}bius transformation with scaling $(2+\sqrt{3})^2$. However, they differ in $D_i$ and $H_i$ strings as described in the text.}
\label{DD2}
\end{figure}

 The eigenvalues of the product of such matrices determine the self-similar scaling\cite{IIS2}.  The golden ( $n=1$)  and the diamond ($n=2)$  hierarchies described above are characterized by $D_3^2$ and $D_3^2D_2$,  encoded in the eigenvalues $ ((\frac{3 \pm \sqrt{5}}{2})^2, 1, 1)$ and $(2\pm \sqrt{3})^2, 1, 1)$ of these string of matrices.
 
As described earlier  the  self-similar properties of the whole gasket  are  described by the M\"{o}bius transformations $ w = f(z)= \frac{ az+b}{cz+d}$  as shown in equation (\ref{Mmap}). These transformations consisting of $ 2 \times 2$ matrices of integers with unit determinant form  the ``Modular group"  $SL(2,Z)$\cite{ModularG}. 
 Modular group is a sub group of $SL(2,C)$ - a group of $2 \times 2$ matrices with unit determinant, where matrix entries are complex numbers.   $SL(2,C)$  describes the recursive structure of the entire
 $\mathcal{IAG}$. 
 Finally, we note that with every self-similar hierarchy, we can associate a M\"{o}bius transformation as well as a string of $D_i$ matrices whose eigenvalues determine the asymptotic scaling. There is  however an important distinction in  these two characterization of the  self-similar hierarchy. As illustrated in Fig. (\ref{DD2}), two very distinct hierarchies, one describing a nested set of circles while the other is a chain of circles, both exhibiting  same scaling, are represented by the same M\"{o}bius transformation. 
 In contrast, within the Apollonian group, these two hierarchies  are respectively described by
$D_3^2D_2 \equiv h_3h_1$ and $D_2D_3^2 \equiv h_1h_3$. In the first case, the nested set of circles exhibit self-similar kaleidoscopic symmetries. In contrast, the $D_2D_3^2 \equiv h_1h_3$ hierarchy does not exhibit this special symmetry. These differences are reflected in
 the eigenvectors  of $h_3h_1$ and $h_1h_3$ which are respectively given by $  \left( \begin{array}{c}   \sin \frac{\pi}{6}  \\  \cos \frac{\pi}{6}  \\  \end{array}\right)$ and $  \left( \begin{array}{c}   \sin \frac{\pi}{12}  \\  \cos \frac{\pi}{12}  \\  \end{array}\right)$.

 \section{ Conclusion}
 
 Apollonian circle packing as described here using various  number theoretical and geometrical properties, relates this fascinating fractal to the  Ford circles, the Farey tree that generates all rationals as well as the the Pythagorean tree for generating primitive Pythagorean triplets.  The entire framework can  be described elegantly using  M\"{o}bius transformations and the group theory involving
 special linear groups $SL(2,C)$ , $SL(2,Z)$ and the Apollonian group. 
 The central result of this paper,  that the self-similar properties of the entire Apollonian gasket are encoded in the Ford circles, can be visualized geometrically due to various kaleidoscopic symmetries of the gasket as shown in Table (\ref{T1}) and further illustrated in Fig. (\ref{K2}). The self-similar recursions associated with the packing are given by  a  M\"{o}bius map
 ({\ref{Mmap}) forming a modular group $SL(2,Z)$, characterized by the matrix $\mathcal{F^*}$ whose eigenvalues describe the asymptotic scalings.
 Intriguingly, a very special family of
irrational numbers are found to be lurking in the self-similar scaling factors. 
 Number theory unites this family, revealing its special elite status as its members, identified with a single integer $n$ 
can be represented by $[ n+1: \overline{1,n}]$. It includes the golden-mean class, that is a group of irrationals with continued fraction representation $[n; \overline{1}]$ as the golden mean is 
$\frac{1+\sqrt{5}}{2} = [1; \overline{1}]$.
It excludes the silver mean class, namely the irrationals with continued fraction expansion as $[n;\overline{2}]$.

As a final comment, we would like to point out that the proof of  the Descartes theorem for the Ford circles stated in the form of the  property of  the homogeneous functions as described in Appendix B sets the stage for an alternative proof of the Descartes theorem that is more elegant than the previously known  proofs of the theorem\cite {J1}. This is because the Descartes configurations for Ford circles are related to an arbitrary Descartes configuration via a conformal transformation.
\appendix

\section{ M\"{o}bius Transformations}

Conformal maps are functions in complex plane that preserve the angles between curves.
There is a specific family of conformal maps  $w = \frac{az+b}{cz+d}$,  known as M\"{o}bius transformation or linear fractional transformation. Such  transformation maps lines and circles to lines and circles.
These maps can be represented  ( denoted by symbol $\dot{=}$ ) with a matrix:
 The map
\begin{equation}
 w = f(z)= \frac{az+b}{cz+d} \,\ \dot{=} \,\ \left( \begin{array}{cc} a & b \\ c & d \\
\end{array} \right) 
 \end{equation}

This identification of M\"{o}bius map with a matrix is useful because if we
consider two conformal maps: $ w_1= f_1(z)$ and $w_2=f_2(z)$ that respectively correspond to the matrices $C_1$ and $C_2$. Then the composition  $f_1f_2(z)$ corresponds to the matrix $C_1\cdot C_2$.

The  constants of Mobius maps $(a,b,c,d)$  can be determined\cite{wiki} in terms of 
two sets of triplets:  $( z_1, z_2, z_3) $ and   their conformal image $ ( w_1, w_2, w_3)$:

\begin{eqnarray}
 a & = & \det { \left( \begin{array}{ccc} z_{1} w_1 & w_1 &1\\ z_{2} w_2 & w_{2} &1\\  z_{3} w_{3} & w_{3} &1  \end{array}\right) } ,\,\
 b  =  \det {  \left( \begin{array}{ccc} z_{1}w_{1}&z_{1}&w_{1}\\z_{2}w_{2}&z_{2}&w_{2}\\z_{3}w_{3}&z_{3}&w_{3}  \end{array}\right) } \nonumber \\
  c & = & \det {  \left( \begin{array}{ccc} z_{1}&w_{1}&1\\z_{2}&w_{2}&1\\z_{3}&w_{3}&1  \end{array}\right)} ,\,\
 d    =   \det {  \left( \begin{array}{ccc} z_{1}w_{1}&z_{1}&1\\z_{2}w_{2}&z_{2}&1\\z_{3}w_{3}&z_{3}&1  \end{array}\right)}
 \label{abcd}
 \end{eqnarray}
 
 In other words, there is a unique  map that connects two distinct set of triplets. Since any three points determine  a configuration of four mutually tangent circles), as illustrated in Fig. (\ref{P3}), any two Descartes configurations are related by a M\"{o}bius map. These relations are consequence of the invariance of ``cross-ratio",
$ R = \frac{(z_4-z_1)(z_2-z_3)}{(z_2-z_1)(z_4-z_3)}$.
 It turns out that the self-similar hierarchies that share a boundary circle correspond to real value of $R$ that remains invariant as we zoom in the recursive structure. 

\section{ Descartes's Theorem for Ford Circles }

Theorem:  If three variables $\lambda_1$, $\lambda_2$ and $\lambda_3$  are related by,
\begin{equation}
\lambda_1 + \lambda_2 + \lambda_3 =0
\end{equation}
Then they also satisfy,
\begin{equation}
 ( \lambda1^2 + \lambda_2^2 + \lambda_3^2)^2 = 2 ( \lambda_1^4 + \lambda_2^4+ \lambda_3^4).
 \label {lam}
 \end{equation}
  Proof:  Define the standard three homogeneous symmetric polynomials on three variables,
  \begin{eqnarray}
 S_1  =  \lambda_1 + \lambda_2 + \lambda_3,\,\,\
 S_2  =  \lambda_1 \lambda_2 + \lambda_2 \lambda_3 + \lambda_3 \lambda_2 ,\,\,\ 
 S_3  =  \lambda_1 \lambda_2 \lambda_3,
 \label{s3}
 \end{eqnarray}
 and some additional ones that will be useful,
  \begin{eqnarray}
 T_2  =  \lambda_1^2 + \lambda_2^2 + \lambda_3^2 ,\,\,\
  T_4   =  \lambda_1^4 + \lambda_2^4 + \lambda_3^4 ,\,\,\
  U  =   \lambda_1^2 \lambda_2^2 + \lambda_2^2 \lambda_3^2 + \lambda_3^2 \lambda_2^2
 \label{U}
 \end{eqnarray}
 By direct multiplication, we have,
  \begin{eqnarray}
 S_1^2   =  T_2 + 2S_2 ,\,\,\
  T_2^2  =  T_4 + 2U ,\,\,\
 S_2^2   =  U + 2 S_1 S_3.
 \label{m3}
 \end{eqnarray}
  We now set $S_1=0$, that gives $T_2 = -2 S_2$. Squaring both sides, and then eliminating U we get,
  \begin{eqnarray}
 T_2 ^2  =  4 S_2^2 
   = 4 U
  = 2 T_2^2 -2T_4
 = 2 T_4
  \label{ed}
  \end{eqnarray}
  $T_2^2= 2T_4$ gives the required equation ( \ref{lam}).  Applying the above formulas for  the Farey relation $ q_c = q_R+q_L$,  we get the Descartes theorem for Ford circles.


\begin{thebibliography}{99}
\bibitem{IAG} D. Mackenzie, Am. Sci. 98 ( 2010 ) 10.

\bibitem{AP} Ronald.Graham,Jeffrey C.Lagarias, Colin L.Mallows, Allan R. Wilks and Catherine H. Yan, Apollonian circle packings: geometry and group theory I. Apollonian group, Discrete and Computational Geometry 34 (2005), 547–585; Ronald L. Graham, Jeffrey C. Lagarias, Colin L. Mallows, Allan R. Wilks and Catherine H. Yan, Apollonian circle packings: number theory, J. Number Theory 100 (2003), 1–45.

\bibitem{Cmap} Ogilvy, C. S. (1990). Excursions in Geometry. Dover. pp. 54 D 55. ISBN 0-486-26530-7;
Also see,  \url{<https://en.wikipedia.org/wiki/Mobius_transformation>}

\bibitem{prime} Peter Sarnak,
 ``Integral Apollonian Packings", 
The American Mathematical Monthly 118(4), March $2011$.

\bibitem{J1} Jerzy Kocik, Proof of Descartes circle formula and its generalization, clarified (arXiv:0706.0372).

\bibitem{J2} Jerzy Kocik, On a Diophantine equation that generates all integral Apollonian gaskets, ISRN
Geometry, 348618 (2012).


\bibitem{Hof}
Hofstadter, D R, 1976,
{\sl Energy-levels and Wave-functions for Bloch electrons in Rational and Irrational Magnetic Fields}, 
{\it Phys.Rev.B}, {\bf 14}, 2239-49.

\bibitem{QHE} Thouless, D J , Kohmoto M, Nightingale M P and den Nijs M, 1982, Quantised Hall Conductance in a Two-Dimensional Periodic Potential, Phys. Rev. Lett., 49, 405-8.

 \bibitem{book} I. I. Satija, {\it Butterfly in the Quantum World} (IOP Concise,
Morgan and Claypool, San Raffael, CA, 2016), Chap. 10.

\bibitem{IIS2} I.I. Satija, {\it J. Phys. A}: Math. Theor. {\bf 54} ( 2021),  025701.
\bibitem{IIS1} I. I. Satija, 2016,
{\sl A tale of two fractals: The Hofstadter butterfly and the integral Apollonian gaskets},
{\it Eur. Phys. J. - Special Topics}, {\bf 225}, 2533-47.

\bibitem{SW} I. Satija and M. Wilkinson,  ``Nests and chains of Hofstadter butterflies",  {\it J. Phys A} ,  {\bf 53},  085703, 2020.

\bibitem{4circle} 
Using $C_1$ as the origin and $C1C2$ as the $x$-axis, 
the coordinates $(x_4,y_4)$ of the center  of the outermost circle of curvature $\kappa_4$ is given by,
$x_4 = \frac{1}{\kappa_1} - \frac{\kappa_2-\kappa_1}{\kappa_4(\kappa_1+\kappa_2)},\,\,\ y_4^2 = (\frac{(\kappa_1-\kappa_4}{\kappa_1 \kappa_4})^2+ x_4^2$.

Also, given the curvatures of the three circles, the sides of the triangles are $( r_1+r_2), (r_1+r_3), (r_2+r_3)$ where $r_i = \frac{1}{\kappa_i}$.

\bibitem{Pchain} Ogilvy, C. S. (1990). Excursions in Geometry. Dover. pp. 54 D 55. ISBN 0-486-26530-7.

\bibitem{Ford}  L. R. Ford, The American Mathematical Monthly, Vol. 45, No. 9 ( 1938) 586-601.
 
\bibitem{RF}This proof is  by Richard Friedberg ( private communication) July 2017.

\bibitem{SatPT}
Satija, I I, 2018, 
{\sl Pythagorean Triplets, Integral Apollonians and The Hofstadter Butterfly,}
arXiv:1802.04585v3 [nlin.CD]


\bibitem{Hall} A. Hall, Geneology of Pythagorean Triads, Math. Gazette 54, No. 390 (1970), 377-379.
\bibitem{J3} Jerzy Kocik, Adv. Appl. Clifford alg. 17 (2007), 793.

\bibitem{ModularG} R. C. Alperin, $PSL2(Z) = Z2*Z3$, Amer. Mathematical Monthly 100 (1993), 385–386.



\end{thebibliography}
\end{document}